\documentclass{article}
\usepackage[leqno]{amsmath}
\usepackage{amsthm,amscd}
\usepackage{amsfonts} 
\usepackage{amssymb} 
\newtheorem{theorem}{Theorem}[section]
\newtheorem{definition}[theorem]{Definition}
\newtheorem{example}[theorem]{Example}

\newtheorem{proposition}[theorem]{Proposition}
\newtheorem{corollary}[theorem]{Corollary}
\newtheorem{remark}[theorem]{Remark}

\newcommand{\A}{{\mathcal A}}

\newcommand{\codim}{\operatorname{codim}}
\newcommand{\Der}{{\rm Der}}
\newcommand{\Poin}{{\rm Poin}}

\newcommand{\nbc}{\mathbf{nbc}}

\begin{document}
\newcommand{\Vanderdet}{\left|
\begin{array}{ccccc}
x_{1} & x_{1}^{3} & . & . &x_{1}^{2l-1}\\
 .    &     .     & . & . &     .      \\
 .    &     .     & . & . &     .      \\
x_{l} & x_{l}^{3} & . & . &x_{l}^{2l-1}     
\end{array}
\right|}

\newcommand{\Vanderdetred}{\left|
\begin{array}{ccccc}
1 & x_{1}^{2} & . & . &x_{1}^{2l-2}\\
 .    &     .     & . & . &     .      \\
 .    &     .     & . & . &     .      \\
1 & x_{l}^{2} & . & . &x_{l}^{2l-2}     
\end{array}
\right|}

\title{Algebras generated by reciprocals of linear forms}
\author{{\sc Hiroaki Terao}
\footnote{partially supported by 
the Grant-in-aid for scientific research (No.1144002), 
the Ministry of Education, Sports, Science and Technology, Japan
}\\
{\small \it Tokyo Metropolitan University, Mathematics Department}\\
{\small \it Minami-Ohsawa, Hachioji, Tokyo 192-0397, Japan}
}
\date{}
\maketitle

\begin{abstract}

Let $\Delta$ be a finite set of nonzero linear forms in several
variables with coefficients in a field $\mathbf K$ of characteristic
zero.
Consider the $\mathbf K$-algebra $C(\Delta)$ of rational functions generated by
$\{1/\alpha \mid \alpha \in \Delta \}$.
Then the ring $\partial(V)$ of differential operators with constant 
coefficients naturally acts on $C(\Delta)$.
We study the graded $\partial(V)$-module structure of $C(\Delta)$.
We especially find standard systems of minimal generators and a
combinatorial formula for the Poincar\'e series of $C(\Delta)$.
Our proofs are based on a theorem by Brion-Vergne \cite{brv1}
and results by Orlik-Terao \cite{ort2}.

\noindent
{\it Mathematics Subject Classification (2000): 32S22, 13D40, 13N10, 52C35} 
\end{abstract}

\bigskip
\setcounter{section}{0}
\setcounter{equation}{0}
\section{Introduction and main results}
~\bigskip

Let $V$ be a vector space of dimension $\ell$ over a field $\mathbf K$
of characteristic
zero.
Let $\Delta$ be a finite subset of the dual space $V^*$ of $V$.
We assume that $\Delta$ does not contain the zero vector
and that no two vectors are proportional throughout
this paper. Let $S=S(V^*)$ be the symmetric algebra of $V^*$.
It is regarded as the algebra of polynomial functions on $V$.
Let $S_{(0)}$ be the field of quotients of $S$, 
which is the field of rational
functions on $V$.

\begin{definition} 
Let $C(\Delta)$ be the $\mathbf K$-subalgebra of $S_{(0)}$ generated
by the set
$$
\{ \frac1{\alpha} \mid \alpha \in \Delta \}.
$$
\end{definition}
Regard $C(\Delta)$ as a graded $\mathbf K$-algebra with 
$\deg (1/\alpha)=1$ for $\alpha \in \Delta$.

\begin{definition} 
Let $\partial(V)$ be the $\mathbf K$-algebra of differential operators
with constant coefficients. 
Agree that the constant multiplications are in
 $\partial(V)$: ${\bf K} \subset \partial(V)$. 
\end{definition}

If $x_1, \cdots, x_\ell$ are a basis for $V^*$,  then $\partial(V)$
is isomorphic to the polynomial algebra $\mathbf {K}[\partial/\partial
x_1, \cdots, \partial/\partial x_\ell$]. Regard $\partial(V)$ as a graded
$\mathbf K$-algebra with $\deg (\partial/\partial x_i)=1 \quad 
(1\leq i\leq\ell)$. 
It naturally acts on $S_{(0)}$.
We regard $C(\Delta)$ as a graded $\partial(V)$-module.
In this paper we study the $\partial(V)$-module structure of $C(\Delta)$.
In particular, we find systems of minimal generators (Theorem ~\ref{1.3})
and a combinatorial formula for the Poincar\'e 
(or Hilbert) series Poin$(C(\Delta), t)$ of $C(\Delta)$
(Theorem ~\ref{1.4}).

In order to present our results we need several definitions. 
Let ${\bf E}_{p}(\Delta)$ be the set of all $p$-tuples composed of
elements of $\Delta$. Let    
${\bf E}(\Delta) := \bigcup_{p\geq 0} {\bf E}_{p}(\Delta)$.
The union is disjoint.
Write $\prod {\cal E} := 
\alpha_{1} \dots \alpha_{p} \in S$ when ${\mathcal E}=
(\alpha_{1}, \dots, \alpha_{p})\in {\bf E}_{p}(\Delta)   $.
Then one can write
\[
C(\Delta) 
=
\sum_{{\mathcal E}\in{\bf E}(\Delta)}
{\mathbf K}\left(\prod {\mathcal E}\right)^{-1}.
\]
Let
\begin{eqnarray*} 
{\bf E}^{i}(\Delta) &=& \{{\cal E}\in {\bf E}(\Delta) \mid {\cal E} 
\text{~is linearly independent}\},\\ 
{\bf E}^{d}(\Delta) &=& \{{\cal E}\in {\bf E}(\Delta) \mid {\cal E} 
\text{~is linearly dependent}\}.
\end{eqnarray*} 
Note that ${\mathcal E}\in {\bf E}^{d}(\Delta)   $ if ${\mathcal E} $ 
contains a repetition.   
In a special lecture at the Japan Mathematical
Society in 1992, K. Aomoto suggested the study of the finite-dimensional
graded $\mathbf K$-vector space
\[
AO(\Delta) :
=
\sum_{{\mathcal E}\in{\bf E}^{i}(\Delta)}
{\mathbf K}\left(\prod {\mathcal E}\right)^{-1}.
\]
Let
$$
\mathcal{A}(\Delta)=\{\ker(\alpha)\mid\alpha \in \Delta\}.
$$
Then $\mathcal{A}(\Delta)$ is a (central) arrangement of hyperplanes
\cite{ort1} in $V$. 
K. Aomoto conjectured, when $\mathbf{K}=\mathbf{R}$, that the dimension of
$AO(\Delta)$ is equal to the number of connected components of
$$
M(\mathcal{A}(\Delta)):=V\setminus \bigcup_{H \in \mathcal{A}(\Delta)}H.
$$
This conjecture was verified  in \cite{ort2}, where
explicit $\mathbf K$-bases for $AO(\Delta)$ were
constructed.
This paper can be considered as a sequel to \cite{ort2}.
(It should be remarked that constructions in \cite{ort2} were generalized
for oriented matroids by R. Cordovil \cite{cor1}.)
We will prove the following

\begin{theorem}
\label{1.3} 
Let $\mathcal B$ be a $\mathbf K$-basis for $AO(\Delta)$.
Let $\partial(V)_+$ denote the maximal ideal of $\partial(V)$ 
generated by the homogeneous elements of
degree one.
Then 

(1)
the set $\mathcal B$ is a system of minimal genarators
for the $\partial(V)$-module $C(\Delta)$,

(2)
$C(\Delta) =\partial(V)_{+}C(\Delta) \oplus AO(\Delta)$, and

(3)
$\partial(V)_{+} C(\Delta)
=
\sum_{{\mathcal E}\in{\bf E}^{d}(\Delta)}
{\mathbf K}\left(\prod {\mathcal E}\right)^{-1}.
$
In paticular, $\partial(V)_{+} C(\Delta)$ is an ideal of $C(\Delta)$.  
\end{theorem}

Let ${\rm Poin}(\mathcal{A}(\Delta), t)$ be the Poincar\'e polynomial
\cite[Definition 2.48]{ort1} of $\mathcal{A}(\Delta)$.
(It is defined combinatorially and is known  to be equal to the the Poincar\'e polynomial
of $M(\mathcal{A}(\Delta))$ when $\mathbf K = \mathbf C$
\cite{ors1} \cite[Theorem 5.93]{ort1}.) 
Then we have

\begin{theorem}
\label{1.4} 
The Poincar\'e series ${\rm Poin}(C(\Delta), t)$ of the graded
module $C(\Delta)$ is equal to 
${\rm Poin}(\mathcal{A}(\Delta), {(1-t)}
^{-1}t)$.
\end{theorem}

In order to prove these theorems we essentially use
a theorem by M. Brion and M. Vergne \cite[Theorem1]{brv1} and
results from \cite{ort2}. 
By Theorem~\ref{1.4} and the factorization theorem
(Theorem \ref{2.4})
in \cite{ter1}, we may easily show the following two corollaries:

\begin{corollary}
\label{1.5} 
If $\mathcal{A}(\Delta)$ is a free arrangement with exponents
$(d_1, \cdots, d_\ell)$ \cite{ter1} \cite[Definitions 4.15, 4.25]{ort1},
then
$$
{\rm Poin}(C(\Delta), t)={(1-t)}^{-\ell} \prod^{\ell}_{i=1}
\{1+(d_{i}-1)t\}.
$$
\end{corollary}

\begin{example}
\label{1.6} 
Let $x_{1},\dots, x_{\ell}   $ be a basis for $V^{*} $.
Let $\Delta = \{x_{i} - x_{j} \mid  1\leq i < j \leq \ell\}$.
Then $\A(\Delta)$ is known to be a free arrangement
with exponents $(0, 1, \dots , \ell-1)$ \cite[Example 4.32]{ort1}.
So, by Corollary~\ref{1.5}, we have
\[
{\rm Poin}(C(\Delta), t) =
(1-t)^{-\ell + 1} (1+t)(1+2t) \dots (1+(\ell - 2)t).  
\]
For example, when $\ell$=3, we have 
\begin{eqnarray*} 
&~&{\rm Poin}\left({\bf K}\left[
\frac{1}{x_{1} - x_{2}},
\frac{1}{x_{2} - x_{3}},
\frac{1}{x_{1} - x_{3}}
\right], t\right) 
=
(1+t)/(1-t)^{2}\\
&=&
1 + 3t + 5 t^{2} + 7 t^{3} + 9 t^4 + \dots,  
\end{eqnarray*} 
which can be easily checked by direct computation.

When  $\A(\Delta)$ is the set of reflecting hyperplanes of 
any (real or complex) reflection group, Corollary
~{1.5} can be applied because 
 $\A(\Delta)$ is known to be a free arrangement \cite{sai1}
\cite{ter2}.
\end{example}

\begin{corollary}
\label{1.6} 
If $\mathcal{A}(\Delta)$ is generic (i, e., $\left|\Delta\right| \geq \ell$ 
and any $\ell$ vectors in
$\Delta$ are linearly independent), then
$$
{\rm Poin}(C(\Delta), t)={(1-t)}^{-\ell} \sum^{\ell-1}_{i=0}
{\binom{|\Delta|-\ell+i-1}{i}}t^i.
$$
\end{corollary}

\section{Proofs} 

In this section we prove Theorems~\ref{1.3} and \ref{1.4}.
For $\varepsilon \in \mathbf{E}(\Delta)$, let $V(\varepsilon)$
denote the set of common zeros of $\varepsilon : 
~V(\varepsilon)=\bigcap^p_{i=1} \ker (\alpha_i)$ when
$\varepsilon=(\alpha_1, \cdots, \alpha_p)$. \\
Define
$$
L=L(\Delta)=\{V(\varepsilon) \mid \varepsilon \in \mathbf{E}
(\Delta)\}.
$$
Agree that $V(\varepsilon)=V$ if $\varepsilon$ is the empty tuple.
Introduce a partial order $\le$ into $L$ by reverse inclusion:
$X\le Y \Leftrightarrow X\supseteq Y$. 
Then $L$ is equal to the intersection lattice of the arrangement
$\mathcal{A}(\Delta)$ \cite[Definition 2.1]{ort1}. 
For $X\in L$, define
$$
\mathbf{E}_X(\Delta):=\{\varepsilon \in \mathbf{E}(\Delta)
\mid V(\varepsilon)=X \}. 
$$
Then
$$
\mathbf{E}(\Delta)=\bigcup_{X \in L}\mathbf{E}_X(\Delta)
~~~\text{(disjoint)}.
$$
Define
$$
C_X(\Delta):=\sum_{\varepsilon \in \mathbf{E}_X(\Delta)} \mathbf{K}
{(\prod \varepsilon)}^{-1}.
$$
Then $C_X(\Delta)$ is a $\partial(V)$-submodule of $C(\Delta)$.
The following theorem is equivalent to Lemma 3.2 in \cite{ort2}.
Our proof is a rephrasing of the proof there.

\begin{proposition}
\label{2.1}
$$
C(\Delta)=\bigoplus_{X\in L} C_X(\Delta).
$$
\end{proposition}
\begin{proof}
It is obvious that $C(\Delta)=\sum_{X \in L}C_X(\Delta)$.
Suppose that $\sum_{X \in L}\phi_X=0$ with $\phi_X \in C_X(\Delta)$.
We will show that $\phi_X=0$ for all $X \in L$.
By taking out the degree $p$ part, we may assume that $\deg \phi_X=p$
for all $X \in L$.
Let $\mathcal{S}=\{X \in L \mid \phi_X \ne 0\}$.
Suppose $\mathcal{S}$ is not empty. Then there exists a minimal element $X_0$
in $\mathcal{S}$ (with respect to the partial order by reverse
inclusion).
Let $X \in \mathcal{S} \setminus \{X_0\}$ and write
$$
\phi_X=\sum_{\varepsilon \in \mathbf{E}_X(\Delta)} c_\varepsilon
{(\prod \varepsilon)}^{-1}
$$
with $c_\varepsilon \in \mathbf{K}$.
Let $\varepsilon \in \mathbf{E}_X(\Delta)$.
Because of the minimality of $X_0$, one has $X_0 \nsubseteq X$.
Thus there exists $\alpha_0 \in \varepsilon$ such that
$X_0 \nsubseteq \ker (\alpha_0)$.
Let $I(X_0)$ be the prime 
ideal of $S$ generated by the polynomial functions vanishing
on $X_0$. Then $\alpha_0 \notin I{(X_0)}$. Thus
$$
(\prod \Delta)^p (\prod \varepsilon)^{-1} \in 
I(X_0)^{p|\Delta_{X_0}|-p+1},
$$
where $\prod \Delta := \prod_{\alpha\in\Delta} \alpha$
and
$\Delta_{X_0}=\Delta \cap I(X_0)$. Multiply
$(\prod \Delta)^p$ to the both sides of
$$
\phi_{X_0}=-\sum_{
\begin{subarray}{l} X \in \mathcal{S} \\
X \ne X_0 \end{subarray}} \phi_X
$$
to get
\begin{eqnarray*}
(\prod \Delta)^p \phi_{X_0} &=& -\sum_{\begin{subarray}{l} 
X \in \mathcal{S} \\
X \ne X_0 \end{subarray}} (\prod \Delta)^p \phi_X \\
&=& -\sum_{\begin{subarray}{l} X \in \mathcal{S} \\ 
X \ne X_0 \end{subarray}} \sum_{\varepsilon \in \mathbf{E}_X(\Delta)}
c_\varepsilon (\prod \Delta)^p (\prod \varepsilon)^{-1} \in
I(X_0)^{p|\Delta_{X_0}|-p+1}.
\end{eqnarray*}
Since $(\prod \Delta)/(\prod \Delta_{X_0}) \in S\setminus I(X_0)$ 
and $I(X_0)^{p|\Delta_{X_0}|-p+1}$ is a primary ideal, one has
$$
(\prod \Delta_{X_0})^p \phi_{X_0} \in I(X_0)^{p|\Delta_{X_0}|-p+1}.
$$
This is a contradiction because
$$
\deg(\prod \Delta_{X_0})^p \phi_{X_0}=p|\Delta_{X_0}|-p.
$$
Therefore $\mathcal{S}=\phi$.
\end{proof}

Next we will study the structure of $C_X(\Delta)$ for each $X \in L$.
Let $AO_X(\Delta)$ be the $\mathbf K$-subspace of $AO(\Delta)$ generated
over $\mathbf K$ by
$$
\{(\prod \varepsilon)^{-1} \mid \varepsilon \in \mathbf{E}(\Delta)^i
\cap \mathbf{E}_X(\Delta) \}.
$$
Then
$$
AO(\Delta)=\bigoplus_{X\in L} AO_X(\Delta)
$$
by Proposition~\ref{2.1}.
Let $\mathcal{B}_X$ be a $\mathbf K$-basis for $AO_X(\Delta)$.
Then we have

\begin{proposition}
\label{2.2}
The $\partial(V)$-module $C_X(\Delta)$ can also be regarded as a free 
$\partial(V/X)$-module with a basis $\mathcal{B}_X$.
In other words, there exists a natural graded isomorphism
$$
\partial(V/X) \otimes_{\mathbf K}  AO_X(\Delta) \simeq C_X(\Delta).
$$
\end{proposition}

\begin{proof}
First assume that $\Delta$ spans $V^*$ and $X=\{\mathbf 0\}$. 
Then $\mathbf{E}(\Delta)^i \cap
\mathbf{E}(\Delta)_X$ is equal to the set of $\mathbf K$-bases
for $V^*$ which are contained in $\Delta$.
Thus $AO_X (\Delta)$ is generated over $\mathbf K$ by
$$
\{(\prod \varepsilon)^{-1} \mid \varepsilon \in \mathbf{E}_\ell (\Delta)
\text{ is a basis for } V\}.
$$
Similarly $C_X$ is spanned over $\mathbf K$ by
$$
\{(\prod \varepsilon)^{-1} \mid \varepsilon \in \mathbf{E}(\Delta)
\text{ spans } V\}.
$$
Then Theorem 1 of \cite{brv1} is exactly the desired result.
Next let $X \in L$ and $\overline{V}=V/X$.
Regard the dual vector space $\overline{V}^*$ as a subspace of $V^*$
and the symmetric algebra $\overline{S}:=S(\overline{V}^*)$ of
$\overline{V}^*$ as a subring of $S$.
Then $\Delta_X:=I(X) \cap \Delta$ is a subset of $\overline{V}^*$ and
$\Delta_X$ spans $\overline{V}^*$.
Consider $AO(\Delta_X)$ and $C(\Delta_X)$ which are both contained
in $\overline{S}_{(0)}$. Note that $C_X(\Delta)$ can be regarded as a
$\partial(V/X)$-module because $\partial(X)$ annihilates
$C_X(\Delta)$.
Denote the zero vector of $\overline{V}$ by $\overline{X}$.
Then it is not difficult to see that
$$
C_{\overline{X}}(\Delta_X) \simeq C_X(\Delta) ~~~(\text{as }\partial
(\overline{V}) \text{-modules}),
$$
$$
AO_{\overline{X}}(\Delta_X) \simeq AO_X(\Delta) ~~~(\text{as }\mathbf{K}
\text{-vector spaces}).
$$
Since there exists a natural graded isomorphism
$$
C_{\overline{X}}(\Delta_X) \simeq \partial(\overline{V}) \otimes_
{\mathbf K} AO_{\overline X}(\Delta_X),
$$
one has 
$$
C_X(\Delta) \simeq \partial(V/X) \otimes_{\mathbf K} AO_X(\Delta).
$$
\end{proof}

\noindent
{\em Proof of Theorem 1.3.}
By Proposition~\ref{2.2}, $C_X(\Delta)$ is generated over
$\partial(V)$ by $AO_X(\Delta)$. Since
$$
C(\Delta)= \bigoplus_{X \in L}C_X(\Delta) ~~~(\text{Proposition 
~\ref{2.1}}),
$$
and
$$
AO(\Delta)= \bigoplus_{X \in L}AO_X(\Delta),
$$
the $\partial (V)$-module
$C(\Delta)$ is generated by $AO(\Delta)$.
So $\mathcal B$ generates
$C(\Delta)$ over $\partial(V)$. 
Define
$$
J(\Delta):= \sum_{\varepsilon \in \mathbf{E}^d (\Delta)}
\mathbf{K} (\prod \varepsilon)^{-1},
$$
which is an ideal of $C(\Delta)$.
Then it is known by \cite[Theorem 4.2]{ort2} that
$$
C(\Delta)=
J(\Delta) \oplus AO(\Delta) ~~~(\text{as } 
\mathbf{K} \text{-vector spaces }).
$$
It is obvious to see that
$$
\partial(V)_+ C(\Delta) \subseteq J(\Delta).
$$
On the other hand, we have
\begin{equation*}
C(\Delta) 
= \partial (V) AO(\Delta)
= \partial(V)_+ AO(\Delta)+ AO(\Delta)
= \partial(V)_+ C(\Delta)+ AO(\Delta).
\end{equation*}
Combining these, we have (2) and (3) at the same time.
By (2), we know that $\mathcal B$ minimally generates
$C(\Delta)$ over $\partial(V)$, which is (1). 
\qed

\bigskip

If $M = \bigoplus_{p\ge 0} M_{p} $ is a graded vector space
with $\dim M_{p} < +\infty \quad (p\ge 0)$, we let
$$\Poin (M, t) = \sum_{p=0}^{\infty}  (\dim M_{p}) t^{p}$$
be its {\bf Poincar\'e} (or {\bf Hilbert) series}. 
Recall \cite[2.42]{ort1} the (one variable) M\"obius function 
$\mu : L(\Delta) \rightarrow {\mathbf Z}$ defined by 
$\mu(V)=1$ and for $X > V$ by $\sum_{Y\leq X}\mu(Y) = 0$.  
Then the {\bf Poincar\'e polynomial} $\Poin (\A(\Delta), t)$ 
of
the arrangement $\A(\Delta)$ is defined by
$$\Poin (\A(\Delta), t)
=
\sum_{X\in L} \mu(X)(-t)^{\codim X}.$$

\begin{proposition}
\label{2.3b} 
{\bf (\cite[Theorem 4.3]{ort2})}
For $X\in L$ we have 
$$\dim AO_{X} (\Delta) = (-1)^{\codim X} \mu(X) 
\text{~and~} \Poin(AO(\Delta), t) = \Poin(\A(\Delta), t).
$$
\end{proposition}

   Recall $C(\Delta)$ is a graded $\partial(V)$-module.
Since $C(\Delta)$ is infinite dimensional, $\Poin(C(\Delta), t)$ 
is a formal power series.  We now prove Theorem~\ref{1.4} 
which gives a combinatorial formula for  $\Poin(C(\Delta), t)$.

\bigskip

\noindent
{\em Proof of Theorem~\ref{1.4}.}
We have 
\[
{\Poin}(C(\Delta), t) = \sum_{X\in L}\Poin(C_{X}(\Delta), t)
=
  \sum_{X\in L}\Poin(\partial(V/X), t) \Poin(AO_{X}(\Delta), t)
\]
by Propositions~\ref{2.1} and \ref{2.2}.  
Since the $\mathbf K$-algebra $\partial(V/X)$ is isomorphic to the polynomial
algebra with $\codim X$ variables, we have
$$
\Poin(C(\Delta), t)=\sum_{X \in L}(1-t)^{-\codim X}
\Poin(AO_X(\Delta), t).
$$
By Proposition~\ref{2.3b}, we have
\[
\Poin(AO_{X}(\Delta), t) = 
(-1)^{\codim X} \mu(X) t^{\codim X}. 
\]
Thus
\begin{eqnarray*}
\Poin(C(\Delta), t) 
&=& 
\sum_{X\in L} 
(-1)^{\codim X} \mu(X) \left(\frac{t}{1-t}\right)^{\codim X}\\
&=&
\Poin(\A(\Delta), (1-t)^{-1} t). \qed
\end{eqnarray*} 

Let $\Der$ be the $S$-module of derivations :
$$
\Der =\{\theta \mid \theta:S \to S \text{ is a } \mathbf{K}\text
{-linear derivations}\}.
$$
Then $\Der$ is naturally isomorphic to $S {\bigotimes}_{\mathbf K}V$.
Define
$$
D(\Delta)=\{\theta \in \Der \mid \theta(\alpha) \in \alpha S 
\text{ for any }\alpha \in \Delta\},
$$
which is naturally an $S$-submodule of $\Der$.
We say that the arrangement $\mathcal{A}(\Delta)$ is {\bf free} if
$D(\Delta)$ is a free $S$-module \cite[Definition 4.15]{ort1}.
An element $\theta\in D(\Delta)$ is said to be {\bf homogeneous of degree
$p$} if $$
 \theta(x)\in S_{p} \text{~for all~}    x \in V^*.
$$
When $\mathcal{A}(\Delta)$ is a free arrangement, let 
$\theta_1, \cdots, \theta_{\ell}$ be a homogeneous basis for $D(\Delta)$.
The $\ell$ nonnegative integers 
$\deg \theta_1, \cdots, \deg \theta_{\ell}$ are called the
{\bf exponents} of $\mathcal{A}(\Delta)$.
Then one has 
\begin{proposition}
{\bf (Factorization Theorem \cite{ter1}, \cite[Theorem 4.137]
{ort1})}
\label{2.4}
If $\mathcal{A}(\Delta)$ is a free arrangement with exponents
$d_1, \cdots, d_{\ell},$ then
$$
\Poin(A(\Delta), t)=\prod^{\ell}_{i=1}(1+d_i t).
$$
\end{proposition}

By Theorem~\ref{1.4} and Proposition~\ref{2.4}, 
we immediately have Corollary~\ref{1.5}. 

\medskip

The arrangement $\mathcal{A}(\Delta)$ is 
{\bf generic}
if $\left|\Delta\right|\geq \ell$ and
any $\ell$
vectors in $\Delta$ are linearly independent. In this case, it is easy
to see that \cite[Lemma 5.122]{ort1}
$$
\Poin(\mathcal{A}(\Delta), t)=(1+t) \sum^{\ell -1}_{i=0}
\binom{|\Delta|-1}{i}t^i.
$$

\noindent
{\em  Proof of Corollary~\ref{1.6}}. 
By Theorem~\ref{1.4}, one has
\begin{eqnarray*}
\Poin(C(\Delta), t) &=& (1+\frac{t}{1-t})\sum^{\ell -1}_{i=0}
\binom{|\Delta|-1}{i} \left({\frac{t}{1-t}}\right)^i \\
&=& (1-t)^{-\ell} \sum^{\ell -1}_{i=0}(1-t)^{\ell -i-1} 
\binom{|\Delta|-1}{i}t^i \\
&=& (1-t)^{-\ell} \sum^{\ell -1}_{i=0} \binom{|\Delta|-1}{i}t^i
\sum^{\ell -i-1}_{j=0} \binom{\ell -i-1}{j}(-1)^jt^j \\
&=& (1-t)^{-\ell} \sum^{\ell -1}_{k=0}t^k \sum^{k}_{j=0}(-1)^j
\binom{|\Delta|-1}{k-j} \binom{\ell -k+j-1}{j}.
\end{eqnarray*}

On the other hand, we have
$$
\sum^{k}_{j=0}(-1)^j \binom{|\Delta|-1}{k-j} \binom{\ell -k+j-1}{j}
= \binom{|\Delta|-\ell +k-1}{k}
$$
by equating the coefficients of $x^k$ in $(1+x)^{|\Delta|-\ell +k-1}$
and $(1+x)^{|\Delta|-1}(1+x)^{-(\ell-k)}$.
This proves the assertion.  \qed

\bigskip

We now consider the nbc (=no broken circuit) 
bases \cite{bjo1} \cite{bjo2} \cite{bjz1} \cite{jat1}
\cite[p.72]{ort2}.
Suppose that $\Delta$ is linearly ordered :
$\Delta=\{\alpha_1, \cdots, \alpha_n\}$. 
Let $X\in L$  with $\codim  X=p$.
Define
\begin{multline*} 
\nbc_X(\Delta)
:=\{\varepsilon \in \mathbf{E}_X (\Delta) \mid
\varepsilon=(\alpha_{i_1}, \cdots, \alpha_{i_p}), i_1 <\cdots<i_p,\\
\text{ contains no broken circuits}\}.
\end{multline*} 
Let
$
{\mathcal B}_X
=\{(\prod \varepsilon)^{-1} \mid \varepsilon \in \nbc_X (\Delta)\}
$
for $X\in L$.
Then we have

\begin{proposition}
{\bf (\cite[Theorem 5.2]{ort2})} 
\label{2.5}
Let $X\in L$. The set ${\mathcal B}_{X}  $ 
is a $\mathbf K$-basis for $AO_X(\Delta)$.
\end{proposition}

Thanks to Propositions \ref{2.1}, \ref{2.2} and \ref{2.5} we
easily have 
\begin{proposition} 
\label{2.6} 
Let
$
{\mathcal B}=\bigcup_{X\in L}{\mathcal B}_X=\{\phi_1, \cdots, \phi_m\}.
$
Write $supp(\phi_i)=X$ if $\phi_i\in {\mathcal B}_X$. Then,
for any $\phi \in C(\Delta)$ and $j \in \{1, \cdots, m\}$,
there uniquely exists $\theta_j \in \partial(V/supp(\phi_j))$ such that
$$
\phi=\sum^m_{j=1} \theta_j (\phi_j).  
$$
\end{proposition} 

\begin{remark} 
Suppose that $\Delta$ spans $V^{*} $ and that
$AO_{\{{\bf 0} \}} (\Delta)
=
\sum_{j=1}^{q} {\bf K} \phi_{j}  $,
where $q = \left|\mu(\{{\bf 0}\})\right|$.
Then the mapping 
$$\phi \mapsto \sum_{j=1}^{q} \theta_{j}^{(0)} (\phi_{j}) 
\in AO_{\{{\bf 0} \}} (\Delta)$$
is the restriction to $C(\Delta)$ of the Jeffrey-Kirwan residue
\cite[Definition 6]{brv1} \cite{sze1}.
Here $\theta_{j}^{(0)}  $ is the degree zero part of
$\theta_{j} $ $(j = 1, \dots , q)$.  
\end{remark}

\end{document}